\documentclass[12pt]{article}
\usepackage{amsthm,amsfonts,amssymb,latexsym}

\theoremstyle{plain}
\newtheorem{corollary}{Corollary}

\newtheorem{theorem}[corollary]{Theorem}
\newtheorem{lemma}[corollary]{Lemma}
\newtheorem*{fact*}{Fact}

\theoremstyle{definition}
\newtheorem*{definition}{Definition}

\theoremstyle{remark}
\newtheorem{remark}{Remark}
\newtheorem{example}{Example}
\newtheorem*{problem*}{Problem}

\def\proofname{Proof}
\def\pf{\ifvmode\else\newline\fi\noindent\textsc{\proofname:\ }}
\def\qed{\mbox{ $\Box$}}
\def\eqref#1{\mbox{\rm (\ref{#1})}}
\renewenvironment{cases}{\left\{\begin{array}{ll}}{\end{array}\right.}
\def\txt#1{\quad\mbox{#1}\quad}
\def\sprod#1#2{\left\langle#1,#2\right\rangle}

\def\E{{\mathbb E}}
\def\F{{\mathbb F}}
\def\N{{\mathbb N}}
\def\R{{\mathbb R}}
\def\sgn{\mathop{\rm sgn}}
\def\Jparam{\mbox{\boldmath$J\kern-0.13em$}}
\def\Sparam{\mbox{\boldmath$S\kern-0.08em$}}
\let\epsilon\varepsilon

\def\decreaseitemsep{\itemsep 0pt plus 2pt minus 1pt}

\begin{document}
\def\footnotemark{}
\title{Superreflexivity and $J$-convexity\\of Banach
  spaces\thanks{Research supported by German Academic Exchange Service
    (DAAD)}}
\author{J\"org Wenzel}
\maketitle

\begin{abstract}
  A Banach space $X$ is superreflexive if each Banach space $Y$ that
  is finitely representable in $X$ is reflexive. Superreflexivity is
  known to be equivalent to $J$-convexity and to the non-existence of
  uniformly bounded factorizations of the summation operators $S_n$
  through $X$.

  We give a quantitative formulation of this equivalence.

  This can in particular be used to find a factorization of $S_n$
  through $X$, given a factorization of $S_N$ through $[L_2,X]$, where
  $N$ is `large' compared to $n$.
\end{abstract}

\section{Introduction}
\label{sec:intro}

Much of the significance of the concept of superreflexivity of a
Banach space $X$ is due to its many equivalent characterizations, see
e.~g. Beauzamy \cite[Part 4]{bea85}.

Some of these characterizations allow a quantification, that makes
also sense in non superreflexive spaces. Here are two examples.

\begin{definition}
  Given $n$ and $0 < \epsilon < 1$, we say that a Banach space $X$
  is \emph{$J(n,\epsilon)$-convex}, if for all elements
  $z_1,\dots,z_n\in U_X$ we have
  \[ \inf_{1\leq k\leq n} \Big\| \sum_{h=1}^k z_h - \sum_{h=k+1}^n z_h
  \Big\| < n(1-\epsilon).
  \]

  We let $\Jparam_n(X)$ denote the infimum of all $\epsilon$, such
  that $X$ is not $J(n,\epsilon)$-convex.
\end{definition}

\begin{definition}
  Given $n$ and $\sigma\geq1$, we say that a Banach space $X$
  \emph{factors the summation operator $S_n$} with norm $\sigma$, if
  there exists a factorization $S_n=B_nA_n$ with $A_n:l_1^n\to X$ and
  $B_n:X\to l_\infty^n$ such that $\|A_n\|\,\|B_n\|=\sigma$.

  We let $\Sparam_n(X)$ denote the infimum of all $\sigma$, such that
  $X$ factors $S_n$ with norm $\sigma$.
\end{definition}

Here, the summation operator $S_n:l_1^n \to l_\infty^n$ is given by
\[ (\xi_k) \mapsto (\sum_{h=1}^k \xi_h)
\]
and $U_X$ denotes the unit ball of the Banach space $X$.

It is known that a Banach space is superreflexive if and only if it is
$J(n,\epsilon)$-convex for some $n$ and $\epsilon > 0$, or
equivalently, if it does not factor the summation operators with
uniformly bounded norm; see James \cite[Th.~5, Lem.~B]{jam67}, and
Sch\"affer/Sundaresan \cite[Th.~2.2.]{schaf70}.

Using the terminology introduced above, this can be reformulated as
follows:

\begin{theorem}
  \label{th:james}
  For a Banach space $X$ the following properties are equivalent:
  \begin{enumerate}
    \decreaseitemsep
  \item $X$ is not superreflexive.
  \item For all $n\in\N$ we have $\Jparam_n(X)=0$.
  \item There is a constant $\sigma\geq1$ such that for all $n\in\N$
    we have $\Sparam_n(X)\leq \sigma$.
  \item For all $n\in\N$ we have $\Sparam_n(X)=1$.
  \end{enumerate}
\end{theorem}

There are two conceptually different methods to prove that $X$ is
superreflexive if and only if $[L_2,X]$ is. The one is to use Enflo's
renorming result \cite[Cor.~3]{enf72}, which is not suited to be
localized, the other is the use of $J$-convexity, see Pisier
\cite[Prop.~1.2]{pis75}. It turns out that for fixed $n$
\begin{equation} \label{eq:j-convex}
  \Jparam_n([L_2,X]) \leq \Jparam_n(X) \leq 4n^2\Jparam_n([L_2,X]).
\end{equation}

Similar results hold also in the case of $B$-convexity; see
\cite[p.~30]{ros76}.
\begin{theorem}
  \label{th:b-convex}
  If for some $n$ and all $\epsilon>0$, $[L_2,X]$ contains
  $(1+\epsilon)$ isomorphic copies of $l_1^n$, then $X$ contains
  $(1+\epsilon)$ isomorphic copies of $l_1^n$.
\end{theorem}
\begin{theorem}
  \label{th:c-convex}
  If for some $n$ and all $\epsilon>0$, $[L_2,X]$ contains
  $(1+\epsilon)$ isomorphic copies of $l_\infty^n$, then $X$ contains
  $(1+\epsilon)$ isomorphic copies of $l_\infty^n$.
\end{theorem}
On the other hand, no result of this kind for the factorization of
$S_n$ is known, i.~e. if for some $n$ and all $\epsilon>0$, $[L_2,X]$
factors $S_n$ with norm $(1+\epsilon)$, does it follow that $X$
factors $S_n$ with norm $(1+\epsilon)$?

Assuming $\Sparam_n([L_2,X])\leq \sigma$ for some constant $\sigma$
and \emph{all} $n\geq1$, one can use Theorem \ref{th:james} to obtain
that $\Jparam_n([L_2,X])=0$ for \emph{all} $n\geq1$ and consequently
$\Sparam_n(X)=1$.

The intent of our paper is to keep $n$ fixed in this reasoning.
Unfortunately, we don't get a result as smooth as Theorems
\ref{th:b-convex} and \ref{th:c-convex}. Instead, we have to consider
two different values $n$ and $N$. If $\Sparam_N([L_2,X])=\sigma$ for
some `large' $N$, then $\Sparam_n(X)\leq(1+\epsilon)$ for some `small'
$n$.  To make this more precise, let us introduce the iterated
exponential (or TOWER) function $P_g(m)$. We let
\[ P_0(m) := m \txt{and} P_{g+1}(m) := 2^{P_g(m)}.
\]
We will prove the following two theorems.

\begin{theorem}
  \label{th:1}
  For fixed $n\in\N$, $\Jparam_n(X)=0$ implies $\Sparam_n(X)=1$.
\end{theorem}
\begin{theorem}
  \label{th:2}
  For fixed $n\in\N$, $\epsilon>0$ and $\sigma\geq1$ there is a number
  $N(\epsilon,n,\sigma)$, such that $\Sparam_N(X)\leq \sigma$ implies
  $\Jparam_n(X)<\epsilon$. The number $N$ can be estimated by
  \[ N\leq P_m(cn),
  \]
  where $m$ and $c$ depend on $\sigma$ and $\epsilon$ only.
\end{theorem}
\begin{remark}\label{r-1}
  The proof of Theorem~\ref{th:1} will in fact show that if
  $\Jparam_n(X)\leq\frac1{n2^{n+1}}$ then
  $\Sparam_n(X)\leq1+2n^2\Jparam_n(X)$. This more quantitative version
  will be needed to prove Corollary~\ref{cor:sumop}.
\end{remark}

Using \eqref{eq:j-convex}, we obtain the following consequence.
\begin{corollary}
  \label{cor:sumop}
  For fixed $n\in\N$, $\epsilon>0$, and $\sigma\geq 1$ there is a
  number $N(\epsilon,n,\sigma)$ such that $\Sparam_N([L_2,X])\leq
  \sigma$ implies $\Sparam_n(X)\leq (1+\epsilon)$. The number $N$ can
  be estimated by
  \[
    N\leq P_m(cn),
  \]
  where $m$ and $c$ depend on $\sigma$ and $\epsilon$ only.
\end{corollary}

The estimate in Theorem \ref{th:2} seems rather crude, and we have no
idea, whether or not it is optimal.

\section{Proofs}
\label{sec:proofs}

First of all, we list some elementary properties of the sequences
$\Sparam_n(X)$ and $\Jparam_n(X)$.

\begin{fact*} \label{fact:1}
\begin{enumerate}
  \decreaseitemsep
\item \label{enum:1}
  The sequence $(\Sparam_n(X))$ is non-decreasing.
\item \label{enum:2}
  $1\leq \Sparam_n(X) \leq (1+\log n)$ for all infinite
  dimensional Banach spaces $X$.
\item \label{enum:3}
  The sequence $(n\Jparam_n(X))$ is non-decreasing.
\item \label{enum:4}
  For all $n,m\in\N$ we have $\Jparam_n(X)\leq \Jparam_{nm}(X) \leq
  \Jparam_n(X)+1/n$.
\item \label{enum:5}
  If $\Jparam_n(X)\to0$ then for all $n\in\N$ we have
  $\Jparam_n(X)=0$.
\item \label{enum:6}
  $\Jparam_n(\R)\geq 1-1/n$ for all $n\in\N$.
\item \label{enum:7}
  If $q$ and $\epsilon$ are related by
  $\epsilon\geq(1-\epsilon)^{q-1}$ then $\Jparam_n(l_q)\leq 4\epsilon$
  for all $n\in\N$.
\end{enumerate}
\end{fact*}
\pf
The monotonicity properties \ref{enum:1} and \ref{enum:3} are trivial.

The bound for $\Sparam_n(X)$ in \ref{enum:2} follows from the fact
that the summation operator $S_n$ factors through $l_2^n$ with norm
$(1+\log n)$ and from Dvoretzky's Theorem.

To see \ref{enum:4} assume that $X$ is $J(n,\epsilon)$-convex. Given
$z_1,\dots,z_{nm}\in U_X$, let
\[ x_h := \frac1m \sum_{k=1}^m z_{(h-1)m+k} \txt{for $h=1,\dots,n$.}
\]
Then
\[ \inf_{1\leq k\leq nm} \Big\| \sum_{h=1}^k z_h - \sum_{h=k+1}^{nm}
z_h \Big\| \leq
m \, \inf_{1\leq k\leq n} \Big\| \sum_{h=1}^k x_h - \sum_{h=k+1}^n x_h
\Big\| < mn(1-\epsilon),
\]
which proves that $X$ is $J(nm,\epsilon)$-convex, and consequently
$\Jparam_n(X)\leq \Jparam_{nm}(X)$.

Assume now that $X$ is $J(nm,\epsilon)$-convex. Given
$z_1,\dots,z_n\in U_X$, let
\[
\begin{array}{ccccccc}
  x_1          & = & \dots & = & x_m    & := & z_1    \\
  \vdots       &   &       &   &        &    & \vdots \\
  x_{(n-1)m+1} & = & \dots & = & x_{nm} & := & z_n.   \\
\end{array}
\]
If
\[ \inf_{1\leq k\leq nm} \Big\| \sum_{h=1}^k x_h - \sum_{h=k+1}^{nm}
x_h \Big\| \txt{is attained for $k_0$,}
\]
there is $l\in\{0,\dots,n\}$ such that $m/2 + (l-1)m<k_0\leq m/2 +
lm$, hence
\[ \Big\| \sum_{h=1}^{k_0} x_h - \sum_{h=k_0+1}^{nm} x_h \Big\| \geq
\Big\| \sum_{h=1}^{lm} x_h - \sum_{h=lm+1}^{nm} x_h \Big\| - 2
\sum_{h\in I} \|x_h\|,
\]
where $I=\{k_0+1,\dots,lm\}$ or $I=\{lm+1,\dots,k_0\}$ according to
whether $k_0\leq lm$ or $k_0>lm$. It follows that
\[ nm(1-\epsilon) > m \inf_{1\leq k\leq n} \Big\| \sum_{h=1}^k z_h -
\sum_{h=k+1}^n z_h \Big\| - m,
\]
and hence $\Jparam_n(X) \geq \epsilon - 1/n$. This proves
\ref{enum:4}.

\ref{enum:5} is a consequence of \ref{enum:4}.

For \ref{enum:6} and \ref{enum:7} see Section \ref{sec:example}.
\qed

Let us now prove Theorem \ref{th:1}.
{\def\proofname{Proof of Theorem \ref{th:1}}
\pf
}
Let $\epsilon>0$, and $z_1,\dots,z_n\in U_X$ be such that
\[ \inf_{1\leq k\leq n} \Big\| \sum_{h=1}^k z_h - \sum_{h=k+1}^n z_h
\Big\| \geq n(1-\epsilon).
\]

By the Hahn-Banach theorem, we find $y_k\in U_{X^\ast}$ such that
\[ n(1-\epsilon) \leq \sum_{h=1}^k \sprod{z_h}{y_k} - \sum_{h=k+1}^n
\sprod{z_h}{y_k}.
\]

Obviously $|\sprod{z_h}{y_k}|\leq1$. If for some $h\leq k$ we even
have
\[ \sprod{z_h}{y_k}<1-n\epsilon,
\] then
\[ n(1-\epsilon) \leq \sum_{l=1}^k \sprod{z_l}{y_k} -
\sum_{l=k+1}^n \sprod{z_l}{y_k} < (n-1) + (1-n\epsilon) =
n(1-\epsilon),
\]
which is a contradiction. Hence
\begin{equation} \label{eq:thm1.1}
  1-n\epsilon \leq \sprod{z_h}{y_k} \leq 1 \txt{for all $h\leq k$.}
\end{equation}
Similarly
\begin{equation} \label{eq:thm1.2}
  1-n\epsilon \leq -\sprod{z_h}{y_k} \leq 1 \txt{for all $h>k$.}
\end{equation}

Let $x_h:=(z_1+z_h)/2$. Then it follows from \eqref{eq:thm1.1} and
\eqref{eq:thm1.2} that there are $x_1,\dots,x_n\in U_X$ and
$y_1,\dots,y_n\in U_{X'}$ so that
\[ \sprod{x_h}{y_k} \in
\begin{cases}
  (1- n\epsilon,1] & \txt{if $h\leq k$,} \\
  (-n \epsilon,+n \epsilon) & \txt{if $h>k$.}
\end{cases}
\]
The assertion now follows from the following distortion lemma.
\qed

\begin{lemma}
  \label{lem:distortion}
  Suppose that for all $\epsilon>0$ there are elements
  $x_1,\dots,x_n\in U_X$ and $y_1,\dots,y_n\in U_{X^\ast}$ such that
  \[ \sprod{x_h}{y_k} \in
  \begin{cases}
    (1-\epsilon,1] & \txt{if $h\leq k$,} \\
    (-\epsilon,+\epsilon) & \txt{if $h>k$.}
  \end{cases}
  \]
  Then $\Sparam_n(X)=1$.
\end{lemma}
\pf
Fix $h\in\{1,\dots,n\}$. Let $\alpha_{lk}:=\sprod{x_l}{y_k}$. Consider
the system of linear equations
\[ \sum_{l=1}^n \alpha_{lk} \xi_l + \alpha_{hk} =
\begin{cases}
  1 & \txt{if $h\leq k$,} \\
  0 & \txt{if $h>k$,}
\end{cases} \qquad k=1,\dots,n
\]
in the $n$ variables $\xi_1,\dots,\xi_n$. Since
$\det(\alpha_{lk})\to1$ if $\epsilon\to0$, this system has a unique
solution $(\xi_l^{(h)})$ with the additional property, that
$|\xi_l^{(h)}|\to0$ as $\epsilon\to0$. Defining $A_n:l_1^n\to X$ by
\[
A_ne_h:= \sum_{l=1}^n x_l\xi_l^{(h)} + x_h,
\]
we get that $\|A_n\|\leq 1+\sup_h \sum_{l=1}^n|\xi_l^{(h)}| \to
1$. Defining $B_n:X\to l_\infty^n$ by
\[
B_nx:=(\sprod{x}{y_k})_{k=1}^n,
\]
we get that $\|B_n\|\leq1$ and $S_n=B_nA_n$. This completes the proof,
since $\Sparam_n(X)\leq\|A_n\|\,\|B_n\|\to1$.  \qed
\begin{remark}
  Note that $\det(\alpha_{lk})\geq 1 - 2^n\epsilon$ and hence
  $|\xi_l^{(h)}|\leq\epsilon/2$ if $\epsilon<1/2^{n+1}$. This shows
  the assertion of Remark~\ref{r-1}.
\end{remark}

\subsection*{Interlude on Ramsey theory}
\label{sec:ramsey}

Our proof of Theorem \ref{th:2} makes massive use of the general form
of Ramsey's Theorem. Therefore, for the convenience of the reader, let
us recall, what it says; see \cite{gra} and \cite{mil86}.

For a set $M$ and a positive integer $k$, let $M^{[k]}$ be the set of
all subsets of $M$ of cardinality $k$.

\begin{theorem}
  Given $r$, $k$ and $n$, there is a number $R_k(n,r)$ such that for
  all $N\geq R_k(n,r)$ the following holds:

  For each function $f:\{1,\dots,N\}^{[k]} \to \{1,\dots,r\}$ there
  exists a subset $M\subseteq\{1,\dots,N\}$ of cardinality at least
  $n$ such that $f(M^{[k]})$ is a singleton.
\end{theorem}

The following estimate for the Ramsey number $R_k(l,r)$ can be found
in \cite[p.~106]{gra}.

\begin{lemma}
  \label{lem:est_ramsey}
  There is a number $c(r,k)$ depending on $r$ and $k$, such that
  \[
  R_k(l,r) \leq P_k(c(r,k)\cdot l).
  \]
\end{lemma}

We can now turn to the proof of Theorem \ref{th:2}.
{
  \def\proofname{Proof of Theorem \ref{th:2}}
  \pf
}
The proof follows the line of James's proof in \cite[Th.~1.1.]{jam64}.
The main new ingredient is the use of Ramsey's Theorem to estimate the
number $N$.

Let $n$, $\epsilon>0$, and $\sigma$ be given. Define $m$ by
\begin{equation} \label{eq:m}
  2m\sigma < \Big( \frac1{1-\epsilon} \Big)^{m-1}
\end{equation}
and let
\begin{equation}\label{eq:N}
  N:= R_{2m}(R_{2m}(2nm+1,m),m),
\end{equation}
where $R$ denotes the Ramsey number introduced in the previous
paragraph.

The required estimate for $N$ then follows from Lemma
\ref{lem:est_ramsey} as follows
\[ N\leq P_{2m}(c_1P_{2m}(c_2 2nm))\leq P_{4m}(c_3n),
\]
where $c_1$, $c_2$, and $c_3$ are constants depending on $m$, which in
turn depends on $\sigma$ and $\epsilon$.

Replacing, e.~g. $\sigma$ by $2\sigma$, we may assume that in fact
$\Sparam_N(X)<\sigma$ in order to avoid using an additional $\delta$
in the notation.  If $\Sparam_N(X)<\sigma$ then there are
$A_N:l_1^N\to X$ and $B_N:X\to l_\infty^N$ such that $S_N=B_NA_N$ and
$\|A_N\|=1$, $\|B_N\|\leq \sigma$. Let $x_h := A_Ne_h$ and $y_k :=
B_N^\ast e_k$.  Note that
\[
\|x_h\| \leq 1, \quad
\|y_k\|\leq \sigma, \txt{and}
\sprod{x_h}{y_k} =
\begin{cases}
  1 & \txt{if $h\leq k$,} \\
  0 & \txt{if $h>k$.}
\end{cases}
\]

For each subset $M\subseteq\{1,\dots,N\}$, we let ${\cal F}_m(M)$
denote the collection of all sequences $\F=(F_1,\dots,F_m)$ of
consecutive intervals of numbers, whose endpoints are in $M$, i.~e.
\[ F_j=\{l_j,l_j+1,\dots ,r_j\}, \quad l_j,r_j\in M, \quad
l_j<r_j<l_{j+1},
\]
for $j=1,\dots,m$. Note that ${\cal F}_m(M)$ can be identified with
$M^{[2m]}$.

The outline of the proof of Theorem~\ref{th:2} is as follows. To each
$\F=(F_1,\dots,F_m)$, we assign an element $x(\F)$ which in fact is a
linear combination of the elements $x_1,\dots,x_N$. Next, we extract a
`large enough' subset $M$ of $\{1,\dots,N\}$, such that all $x(\F)$
with $\F\in{\cal F}_m(M)$ have about equal norm. Finally, we look at
special sequences $\F^{(1)},\dots,\F^{(n)}$ and
$\E^{(1)},\dots,\E^{(n)}$ in ${\cal F}_m(M)$ such that
\[ \Big\| \sum_{h=1}^k x(\F^{(h)}) - \sum_{h=k+1}^n
x(\F^{(h)})\Big\| \geq n \|x(\E^{(k)})\|.
\]
Since $\|x(\E^{(k)})\|\asymp\|x(\F^{(h)})\|$, normalizing the elements
$x(\F^{(h)})$ yields the required elements $z_1,\dots,z_n$ to prove
that $\Jparam_n(X)<\epsilon$.

Let us start by choosing the elements $x(\F)$. For a sequence $\F\in
{\cal F}_m(M)$, we define
\[ S(\F) := \left\{ x=\sum_{h=1}^N \xi_h x_h \colon
  \sup_h|\xi_h|\leq 2, \
  \sprod x{y_l} = (-1)^j \
  \begin{array}{l}
    \mbox{for all $l\in F_j$} \\
    \mbox{and $j=1,\dots,m$}
  \end{array}
\right\}.
\]
By compactness, there is $x(\F)\in S(\F)$ such that
\[ \|x(\F)\| = \inf_{x\in S(\F)} \|x\|.
\]

\begin{lemma} \label{lem:size}
  We have $1/\sigma \leq \|x(\F)\| \leq 2m$ for all $\F\in{\cal
    F}_m(\{1,\dots,N\})$.
\end{lemma}
\pf
Write $F_j=\{l_j,\dots,r_j\}$ and let
\[ x:= -x_{l_1} + 2 \sum_{i=2}^m (-1)^i x_{l_i}.
\]
Then for $l\in F_j$, we have
\[
  \sprod x{y_l}
  =
  -1 + 2 \sum_{i=2}^j (-1)^i \cdot 1 + 2 \sum_{i=j+1}^m (-1)^i \cdot 0 =
  (-1)^j,
\]
hence $x\in S(\F)$ and $\|x(\F)\| \leq \|x\| \leq 2m-1$.

On the other hand,
\[ 1 = |\sprod {x(\F)}{y_{l_1}}| \leq \sigma\|x(\F)\|.
\] Hence $1/\sigma \leq\|x(\F)\|$. \qed

By \eqref{eq:m}, we can write the interval $[1/\sigma,2m]$ as a
disjoint union as follows
\[ \Big[\frac1\sigma,2m \Big]\subseteq \bigcup_{i=1}^{m-1} A_i,
\txt{where} A_i := \frac1\sigma \Bigg[
\Big(\frac1{1-\epsilon}\Big)^{i-1},
\Big(\frac1{1-\epsilon}\Big)^i \Bigg).
\]

For $\F=(F_1,\dots,F_m)\in{\cal F}_m(\{1,\dots,N\})$ and $1\leq j\leq
m$, let
\[ P_j(\F):=(F_1,\dots,F_j)\in{\cal F}_j(\{1,\dots,N\}).
\]
Obviously
\[ \|x(P_{j-1}(\F))\| \leq \|x(P_j(\F))\| \leq 2m \txt{for
  $j=2,\dots,m$.}
\]
It follows that for each $\F\in{\cal F}_m(\{1,\dots,N\})$ there is at
least one index $j$ for which the two values $\|x(P_{j-1}(\F))\|$ and
$\|x(P_j(\F))\|$ belong to the same interval $A_i$. Letting $f(\F)$ be
the least such value $j$, defines a function
\[ f:\{1,\dots,N\}^{[2m]}\to \{1,\dots,m\}.
\]
Applying Ramsey's Theorem to that function, yields the existence of a
number $j_0$ and a subset $L$ of $\{1,\dots,N\}$ of cardinality
$|L|\geq R_{2m}(2nm+1,m)$ such that for all $\F\in{\cal F}_m(L)$ the
two values $\|x(P_{j_0-1}(\F))\|$ and $\|x(P_{j_0}(\F))\|$ belong to
the same of the intervals $A_i$.

Next, for each $\F\in{\cal F}_m(L)$ there is a unique number $i$ for
which the value $\|x(P_{j_0}(\F))\|$ belongs to the interval $A_i$.
Letting $g(\F)$ be that number $i$, defines a function
\[ g:L^{[2m]}\to\{1,\dots,m\}.
\]
Applying Ramsey's Theorem to that function, yields the existence of a
number $i_0$ and a subset $M$ of $L$ of cardinality $|M|\geq 2nm+1$
such that for all $\F\in{\cal F}_m(M)$ we have
\begin{equation} \label{eq:x(F)}
\|x(P_{j_0}(\F))\|\in A_{i_0},
\end{equation}
and hence, by the choice of $j_0$ and $L$, also
\begin{equation} \label{eq:x(E)}
\|x(P_{j_0-1}(\F))\|\in A_{i_0}.
\end{equation}

We now define sequences
\[
\F^{(h)}:=(F_1^{(h)},\dots,F_m^{(h)}) \txt{and}
\E^{(k)}:=(E_1^{(k)},\dots,E_{m-1}^{(k)})
\]
of nicely overlapping intervals.

Write $M=\{p_1,\dots,p_{2nm+1}\}$, where $p_1<p_2<\dots<p_{2nm+1}$ and
define
\[ \F^{(h)} := (F_1^{(h)},\dots,F_m^{(h)}) \in {\cal F}_m(M)
\qquad h=1,\dots,n
\]
as follows
\[
F_j^{(h)}:=
\begin{cases}
    \{p_h,\dots,p_{n+2h-1}\} &
                 \txt{if $j=1$,} \\[2pt]
    \{p_{n(2j-3)+2h},\dots,p_{n(2j-1)+2h-1}\} &
                 \txt{if $j=2,\dots,m-1$,} \\[2pt]
    \{p_{n(2m-3)+2h},\dots,p_{n(2m-1)+h}\} &
                 \txt{if $j=m$.}
\end{cases}
\]
It turns out that
\begin{equation} \label{eq:def_E}
  E_j^{(k)} := \bigcap_{h=1}^k F_{j+1}^{(h)} \cap \bigcap_{h=k+1}^n
  F_j^{(h)} \qquad k=1,\dots,n
\end{equation}
is given by
\[
E_j^{(k)}:=
   \{p_{n(2j-1)+2k},\dots,p_{n(2j-1)+2k+1}\} \txt{if $j=1,\dots,m-1$.}
\]
Hence $(E_1^{(k)},\dots,E_{m-1}^{(k)})\in{\cal F}_{m-1}(M)$. In order
to obtain an element of $\mathcal{F}_m(M)$ we add the auxiliary set
$E_m^{(k)}:=\{p_{2nm},\dots,p_{2nm+1}\}$, this can be done for
$n\geq2$, which is the only interesting case anyway since
$\Jparam_1(X)=0$ for any Banach space $X$. We have
$\E^{(k)}:=(E_1^{(k)},\dots,E_m^{(k)})\in{\cal F}_m(M)$.

The following picture shows the sets $F_j^{(h)}$ and $E_j^{(k)}$ in
the case $n=3$ and $m=4$:

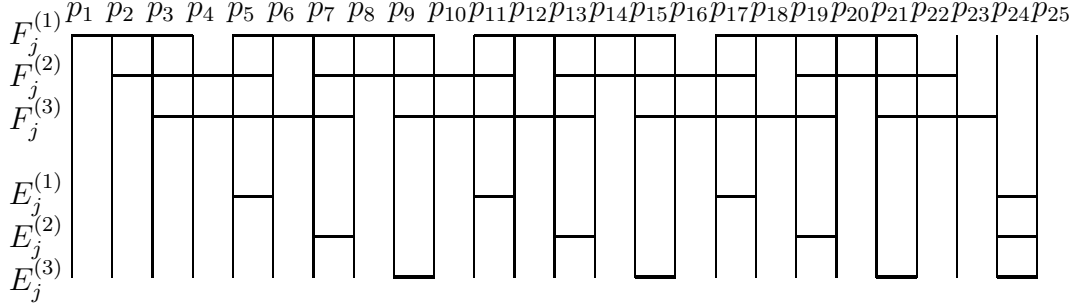
\begin{center}
\setlength{\unitlength}{0.0007in}%

\begin{picture}(4523,2220)(1100,-2220)

\thicklines
\multiput(1200,-1500)(1800,0){3}{\line(1,0){300}}
\multiput(1800,-1800)(1800,0){3}{\line(1,0){300}}
\multiput(2400,-2100)(1800,0){3}{\line(1,0){300}}
\put(6900,-1500){\line(1,0){300}}
\put(6900,-1800){\line(1,0){300}}
\put(6900,-2100){\line(1,0){300}}

\put(0,-300){\line(1,0){900}}
\put(300,-600){\line(1,0){1200}}
\put(600,-900){\line(1,0){1500}}
\multiput(1200,-300)(600,-300){3}{\line(1,0){1500}}
\multiput(3000,-300)(600,-300){3}{\line(1,0){1500}}
\put(4800,-300){\line(1,0){1500}}
\put(5400,-600){\line(1,0){1200}}
\put(6000,-900){\line(1,0){900}}

\thinlines
\multiput(0,-300)(300,0){25}{\line(0,-1){1800}}

\newcounter{interval}

\setcounter{interval}{1}
\multiput(-470,-450)(0,-320){3}{\makebox(0,0)[lb]{%
    $F_j^{(\arabic{interval}\stepcounter{interval})}$}}

\setcounter{interval}{1}
\multiput(-470,-1650)(0,-320){3}{\makebox(0,0)[lb]{%
    $E_j^{(\arabic{interval}\stepcounter{interval})}$}}

\setcounter{interval}{1}
\multiput(-40,-200)(300,0){25}{\makebox(0,0)[lb]{%
    $p_{\arabic{interval}\stepcounter{interval}}$}}

\end{picture}
\end{center}

It follows from \eqref{eq:def_E} that for $1\leq k\leq n$
\[ \frac1n \Big( -\sum_{h=1}^k x(P_{j_0}(\F^{(h)})) + \sum_{h=k+1}^n
x(P_{j_0}(\F^{(h)})) \Big) \in S(P_{j_0-1}(\E^{(k)}))
\]
hence
\[
\Big\| \sum_{h=1}^k x(P_{j_0}(\F^{(h)})) - \sum_{h=k+1}^n
x(P_{j_0}(\F^{(h)})) \Big\| \geq n \|x(P_{j_0-1}(\E^{(k)}))\|.
\]

Let $z_h := \sigma\, (1-\epsilon)^{i_0} \, x(P_{j_0}(\F^{(h)}))$. Then
\[ \Big\| \sum_{h=1}^k z_h - \sum_{h=k+1}^n z_h \Big\| \geq
n\,\sigma \, (1-\epsilon)^{i_0} \, \|x(P_{j_0-1}(\E^{(k)}))\|.
\]
By \eqref{eq:x(F)} we have $\|x(P_{j_0}(\F^{(h)}))\|\in A_{i_0}$,
which implies $\|z_h\|\leq1$. On the other hand, by \eqref{eq:x(E)} we
have $\|x(P_{j_0-1}(\E^{(k)}))\|\in A_{i_0}$, which implies
\[ \Big\| \sum_{h=1}^k z_h - \sum_{h=k+1}^n z_h \Big\| \geq
n\,\sigma \, (1-\epsilon)^{i_0} \, \frac1\sigma
\left( \frac1{1-\epsilon} \right)^{i_0-1} = n\, (1-\epsilon).
\]
Consequently $\Jparam_n(X)\leq \epsilon$. \qed

\section{Problems and Examples}
\label{sec:example}

\begin{example}
  $\Jparam_n(\R)\geq 1-1/n$.
\end{example}
\pf
Let $|\xi_h|\leq1$ for $h=1,\dots,n$. For $k=1,\dots,n$ define
\[ \eta_k:=\sum_{h=1}^k \xi_h - \sum_{h=k+1}^n \xi_h
\] and let $\eta_0:=-\eta_n$. Obviously $|\eta_k-\eta_{k+1}|\leq2$ for
$k=0,\dots,n-1$. Since $\eta_0=-\eta_n$ there exists at least one
$k_0$ such that $\sgn\eta_{k_0} \not= \sgn\eta_{k_0+1}$. Assume that
$|\eta_{k_0}|>1$ and $|\eta_{k_0+1}|>1$, then
$|\eta_{k_0}-\eta_{k_0+1}|>2$, a contradiction. Hence there is $k$
such that $|\eta_k|\leq1$. This proves that
\[ \inf_{1\leq k\leq n} \Big| \sum_{h=1}^k \xi_h - \sum_{h=k+1}^n
\xi_h \Big| \leq 1 = n\frac1n,
\] and hence $\Jparam_n(\R)\geq 1-\frac1n$. \qed

\begin{example}
  If $q$ and $\epsilon$ are related by
  \[ \epsilon \geq (1-\epsilon)^{q-1} \]
  then $\Jparam_n(l_q)\leq 4\epsilon$ for all $n\in\N$.
\end{example}
\pf
Given $\epsilon>0$ find $n_0$ such that
\[ \frac1{n_0} < \epsilon \leq \frac1{n_0-1},
\]
then
\[ \Big( \frac1{n_0} \Big)^{1/q} \geq \Big(1- \frac1{n_0}
\Big)^{1/q} \epsilon^{1/q} \geq 1-\epsilon.
\]
If $n\leq n_0$, choosing
\[ x_h:=( \overbrace{-1,\dots,-1}^h, \overbrace{+1,\dots,+1}^{n-h},
0, \dots),
\]
we obtain
\[ \Big\| \sum_{h=1}^k x_h - \sum_{h=k+1}^n x_h \Big\|_q \geq
\Big\| \sum_{h=1}^k x_h - \sum_{h=k+1}^n x_h \Big\|_\infty = n.
\]
And since
\[ \|x_h\|_q=n^{1/q} \leq n_0^{1/q} \leq 1/(1-\epsilon)
\]
it follows that $\Jparam_n(l_q) \leq \epsilon$.

If $n>n_0$, there is $m\geq2$ such that $(m-1)n_0 < n \leq mn_0$.
Hence, by Properties \ref{enum:3} and \ref{enum:4} in the fact in
Section \ref{sec:proofs} it follows that
\[ \Jparam_n(X) \leq \frac{mn_0}n \Jparam_{mn_0}(X) \leq \frac{mn_0}n
(\Jparam_{n_0} + \frac1{n_0}) \leq \frac{mn_0}n 2\epsilon \leq
4\epsilon. \qed
\]

The main open problem of this article is the optimality of the
estimate for $N$ in Theorem \ref{th:2}.
\begin{problem*}
  Are there $\sigma\geq1$ and $\epsilon>0$ and a sequence of Banach
  spaces $(X_n)$ such that
  \[ \Sparam_{f(n)}(X_n) \leq \sigma \txt{and} \Jparam_n(X_n) \geq
  \epsilon,
  \]
  where $f(n)$ is any function such that $f(n)>n$?

  In particular $f(n) > P_m(n)$, where $m$ is given by \eqref{eq:m}
  would show that the estimate in Theorem \ref{th:2} for $N$ is sharp
  in an asymptotic sense.
\end{problem*}

\vspace{1ex}
\noindent
{\sc Mathematisches Institut, FSU Jena, 07740 Jena, Germany}\\
{\it E--mail:} {\tt wenzel@minet.uni-jena.de}

\end{document}